\pdfoutput=1
\documentclass[12pt,a4paper,reqno]{amsart}
\usepackage{amssymb}

\usepackage{amscd}
\usepackage{enumerate}
\usepackage{graphicx}
\usepackage{siunitx}
\usepackage{tikz-cd}
\usepackage{color}
\usetikzlibrary{arrows}

\usepackage{mdframed}   % for framing
\usepackage{mathabx}

\usepackage{mathtools}%                  http://www.ctan.org/pkg/mathtools
\usepackage[tableposition=top]{caption}% http://www.ctan.org/pkg/caption
\usepackage{booktabs,dcolumn}%           http://www.ctan.org/pkg/dcolumn + http://www.ctan.org/pkg/booktabs

\DeclareCaptionLabelSeparator{separation}{:\quad}
\DeclareFontFamily{OT1}{rsfs}{}
\DeclareFontShape{OT1}{rsfs}{n}{it}{<-> rsfs10}{}
\DeclareMathAlphabet{\mathscr}{OT1}{rsfs}{n}{it}

\theoremstyle{plain}

\newtheorem{theorem}{Theorem}[section]

\theoremstyle{definition}

\newtheorem{definition}[theorem]{Definition}

\usepackage{hyperref}

\usepackage{parskip}
\DeclareRobustCommand{\SkipTocEntry}[5]{}

\usepackage[most]{tcolorbox}
\tcbuselibrary{theorems}
\usepackage{kantlipsum}

\newtcbtheorem{Question}{Question}%
{breakable,colframe=red!50!white!50!black,fonttitle=\bfseries}{}

\begin{document}

%\title[Questions on the Structure of Perfect Matchings]{Questions on the Structure of Perfect Matchings \\ \small{inspired by Quantum Physics Experiments}}

\title[Questions on the Structure of Perfect Matchings in Graphs]{%
  Questions on the Structure of Perfect Matchings
  inspired by Quantum Physics%
}

\author{Mario Krenn}
\address{Vienna Center for Quantum Science \& Technology (VCQ), Faculty of Physics, University of Vienna, Boltzmanngasse 5, 1090 Vienna, Austria.}
\address{Institute for Quantum Optics and Quantum Information (IQOQI), Austrian Academy of Sciences, Boltzmanngasse 3, 1090 Vienna, Austria.}
\address{present address: Department of Chemistry, University of Toronto, Toronto, Ontario M5S 3H6, Canada.}
\address{present address: Vector Institute for Artificial Intelligence, Toronto, Canada.}
\email{mario.krenn@univie.ac.at}

\author{Xuemei Gu}
\address{State Key Laboratory for Novel Software Technology, Nanjing University, 163 Xianlin Avenue, Qixia District, 210023, Nanjing City, China.}
\address{Institute for Quantum Optics and Quantum Information (IQOQI), Austrian Academy of Sciences, Boltzmanngasse 3, 1090 Vienna, Austria.}
\email{xmgu@smail.nju.edu.cn}

\author{D\'{a}niel Solt\'{e}sz}
\address{Alfr\'{e}d R\'{e}nyi Institute of Mathematics, Hungarian Academy of Sciences, 13-15 Re\'{a}ltanoda Street, 1053 Budapest, Hungary.}
\email{solteszd@renyi.hu}


\begin{abstract} We state a number of related questions on the structure of perfect matchings. Those questions are inspired by and directly connected to Quantum Physics. In particular, they concern the constructability of general quantum states using modern photonic technology. For that we introduce a new concept, denoted as \textit{inherited vertex coloring}. It is a vertex coloring for every perfect matching. The colors are inherited from the color of the incident edge for each perfect matching. First, we formulate the concepts and questions in pure graph-theoretical language, and finally we explain the physical context of every mathematical object that we use. Importantly, every progress towards answering these questions can directly be translated into new understanding in quantum physics.
\end{abstract}

\maketitle
%\today

\setcounter{tocdepth}{1}
\tableofcontents
% {\small{\listoftables
% \listoffigures}}

%%%%%%%%%%%%%%%%%%%%%%%%%%%%%%%%%%%%%%%%%%%%%%%%%

\section{Motivation}\label{sec:intro}
A bridge between quantum physics and graph theory has been uncovered recently \cite{krenn2017quantum,gu2018quantum,gu2018quantum3}. It allows to translate questions from quantum physics -- in particular about photonic quantum physical experiments -- into a purely graph theoretical language. The question can then be analysed using tools from graph theory and the results can be translated back and interpreted in terms of quantum physics. The purpose of this manuscript is to collect and formulate a large class of questions that concern the generation of pure quantum states with photons with modern technology. This will hopefully allow and motivate experts in the field to think about these issues.

More concrete, the problems that we present here are concerned with the design of quantum experiments for producing high-dimensional and multipartite entangled quantum states using state-of-the-art photonic technology \cite{pan2012multiphoton}. We start by asking for the generation of Greenberger-Horne-Zeilinger (GHZ) states \cite{greenberger1989going}, and their high-dimensional generalisations \cite{ryu2013greenberger, ryu2014multisetting, lawrence2014rotational, erhard2018experimental}, and further generalize the questions to cover arbitrary pure quantum states.

The paper is organized as follows. In Section 2, we rigorously define the graph theoretic questions that turn out to be relevant in quantum physics. In Section 3 we discuss the correspondence between the all mathematical objects used in Section 2 and quantum experiments.  

\section{Concepts and Questions}\label{sec:concepts}
The type of quantum experiments, that we are interested in, correspond to so-called bi-colored graphs, that are defined as follows. 

\begin{definition}[Edge bi-colored weighted graph] Let $C=\{c_1,\ldots,c_d\}$ be the set of $d\geq 2$ distinct colors. An edge bi-colored weighted graph G=(V(G),E(G)), on $n$ vertices with $d\geq 2$ colors is an undirected, loopless graph where there is a fixed ordering of the vertices $v_1,\ldots,v_n\in V(G)$ and to each edge $e \in E(G)$ a complex weight $w_e$ and an ordered pair of (not necessarily different) colors from $C$ is associated. We say that an edge is monochromatic if two associated colors are not different, otherwise the edge is bi-chromatic. Moreover, if $e$ is an edge incident to the vertices $v_i,v_j \in V(G)$ with $i<j$ and the associated ordered pair of colors to $e$ is $(c_1(e),c_2(e))$ then we say that $e$ is colored $c_1$ at at the endpoint $v_i$ and $c_2$ at the endpoint $v_j$.
\end{definition}

For simplicity, for the rest of the manuscript we appreviate \textit{edge bi-colored weighted graph} by \textit{bi-colored graph}.

The unusual property of bi-colored graphs (compared to other edge-colorings in graph theory) is that edges are allowed to have different colors at different endpoints. The next definition will establish a connection between perfect matchings and vertex colorings of a bi-colored graph.

\begin{definition}[Inherited Vertex Coloring]Let $G$ be a bi-colored graph and let $PM$ denote a perfect matching in $G$. We associate a coloring of the vertices of G with $PM$ in the natural way: for every vertex $v_i$ there is a single edge $e(v_i) \in PM$ that is incident to $v_i$, let the color of $v_i$ be the color of $e(v_i)$ at $v_i$. We call this coloring the inherited vertex coloring (IVC) of the perfect matching $PM$ and denote it by $c$. When all vertices in IVC are colored with only one color, we call $c$ a monochromatic coloring.
\end{definition}

Now we are ready to define how constructive and destructive interference during an experiment is governed by perfect matchings of a bi-colored graph. 

\begin{definition}[Weight of Vertex Coloring] Let $G$ be a bi-colored graph. Let $\mathcal{M}$ be the set of perfect matchings of $G$ which have the coloring $c$ as their inherited vertex coloring. We define the weight of $c$ as 
$$w(c) := \sum_{PM \in \mathcal{M}} \prod_{e \in PM}w_e. $$
Moreover, if $w(c)$=1 we say that the coloring gets unit weight, and if $w(c)$=0 we say that the coloring cancels out. 
\end{definition}

\begin{figure}[!ht]
\centering
\includegraphics[width=\textwidth]{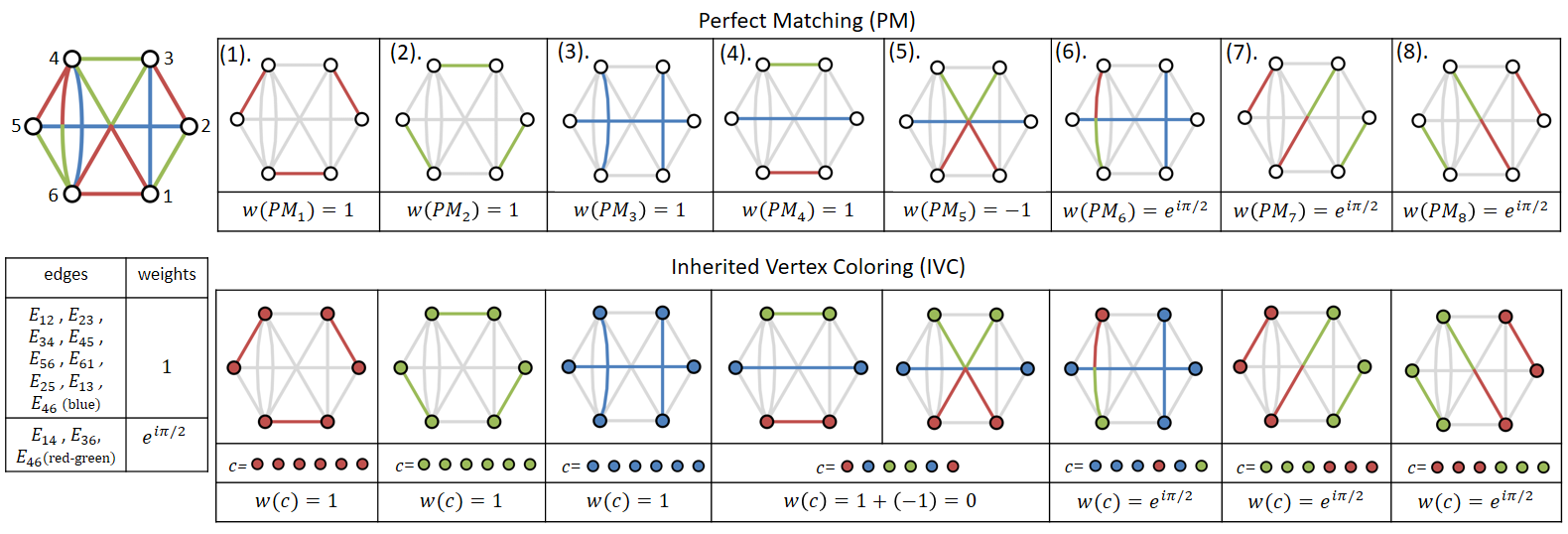}
\caption{Example for inherited vertex coloring and coloring weight. A bi-chromatic weighted edge with one double edge between vertex 4 and 6 is shown on the top left, the edge weights $E_{ij}$ are shown below. On the right top, its eight perfect matchings are shown, and $w(PM_i)$ denotes the product of the edge weights of the perfect matching $PM_i$. The perfect matching 4 and 5 have the same inherited vertex coloring. As $w(c)=w(PM_4)+w(PM_5)=0$, we say this coloring cancels out. There are six remaining IVCs with nonzero weights.}
\label{fig:fig0MultiGraphVertex6}
\end{figure}

An example for a bi-colored graph where some colorings of the vertices get unit weight and some other colorings cancel out can be seen in Figure \ref{fig:fig0MultiGraphVertex6}.

\begin{Question}{monochromatic graph}{} For which values of $n$ and $d$ are there bi-colored graphs on $n$ vertices and $d$ different colors with the property that all the $d$ monochromatic colorings have unit weight, and every other coloring cancels out? We call such a graph \textit{monochromatic}.\label{Question1}
\end{Question}
\begin{figure}[b]
\centering
\includegraphics[width=0.5\textwidth]{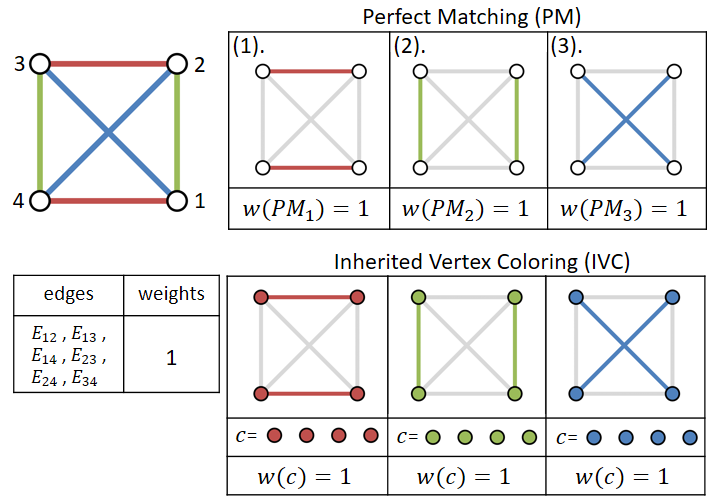}
\caption{A bi-coloring and weight assignment of the edges of $K_4$ that demonstrated that the answer to Question 1 is affirmative for $d=3,n=4$.}
\label{fig:fig1GHZK4}
\end{figure}
The only known values of $n$ and $d$, for which the answer for Question 1 is affirmative, are $d=2$ and $n$ arbitrary even, and $d=3$ ,$n=4$. For $d=2$ and $n$ even an alternately colored (all edges are monochromatic) even cycle $C_n$ suffices with all edge weights being one. For $d=3$, $n=4$ a suitable bi-colored graph can be constructed as follows. Decompose the edges of the complete graph $K_4$ into three disjoint perfect matchings, and let the edges of these matchings be monochromatic, and colored with different color, finally assign weight $w_e=$1 to each edge. It is easy to check that the resulting graph satisfies the conditions of Question 1, see Figure \ref{fig:fig1GHZK4}. Observe that in all known cases we can use weight 1 for each edge. It was shown by Ilya Bogdanov that no other examples are possible with the restriction that all edge weights are positive \cite{bogdanov267013}. The graph in Figure \ref{fig:fig0MultiGraphVertex6} is not monochromatic.

In quantum experiments, one can use additional heralding photons in order to produce a certain state. This concept can be formulated in the following way.

\begin{definition}[$k$-monochromatic colorings] A coloring $c$ is called $k$-monochromatic, if the first $k \leq |V|$ vertices have the same color, and all other vertices are colored (without loss of generality) red.
\end{definition}

\begin{figure}[!ht]
\centering
\includegraphics[width=\textwidth]{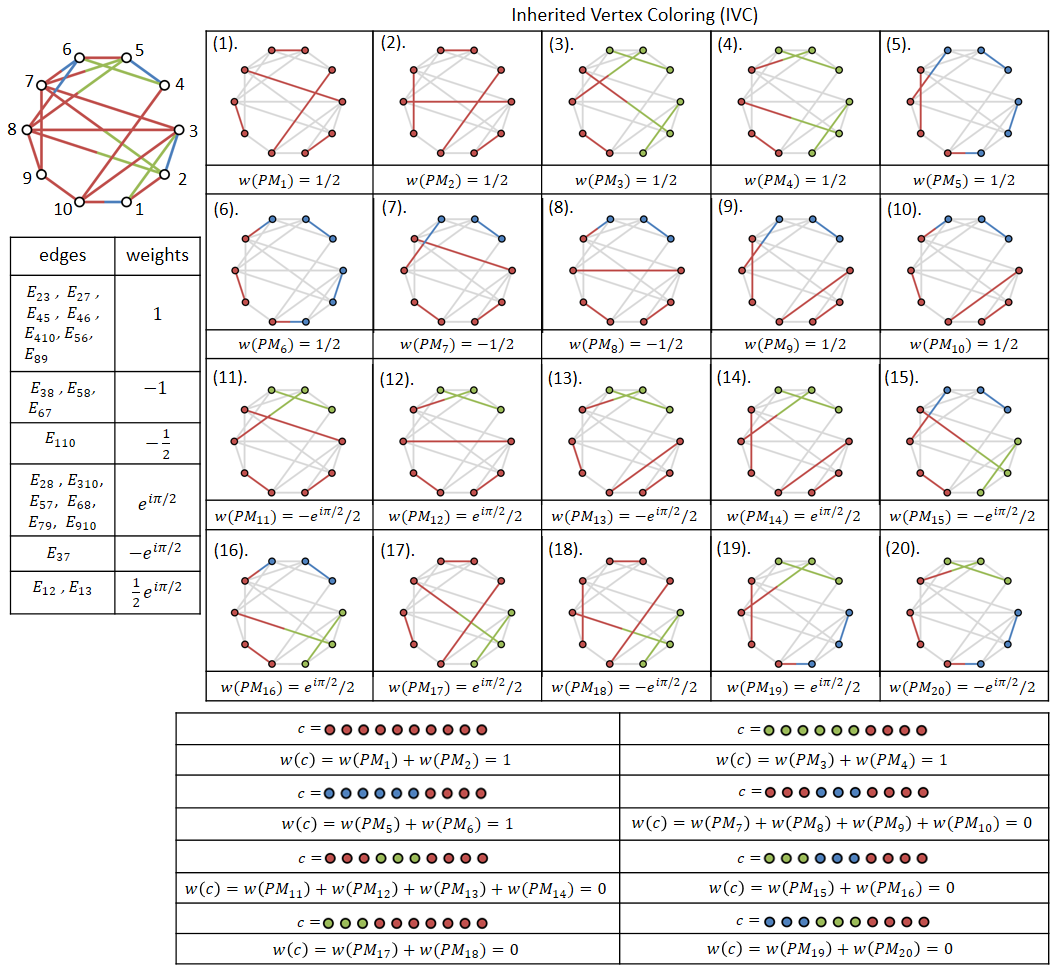}
\caption{The \textit{Erhard graph} is $6-monochromatic$. It is the only known example for $k>4$ and $d\geq3$ satisfying Question 2.}
\label{fig:fig3GHZ4Triggers}
\end{figure}

\begin{Question}{$k$-monochromatic Graph}{} For which values of $n$, $d$ and $k$ are there bi-colored graphs on $n$ vertices and $d$ different colors with the property that all the $d$ $k$-monochromatic colorings have unit weight, and every other coloring cancels out? We call such a graph \textit{k-monochromatic}.\label{Question2}
\end{Question}
The only known example of a $k$-monochromatic graph with $k>4$ and $d\geq3$ is shown in Figure \ref{fig:fig3GHZ4Triggers}. There are three 6-monochromatic colorings, where each has $w(c)=1$. All other colorings are non-6-monochromatic, and have a weight of $w(c)=0$. We call this graph \textit{Erhard graph}\footnote{It is named after Manuel Erhard, who discovered the quantum mechanical technique which has inspired the construction of this graph.}. Note that increasing the number $n$ while keeping $k$ constant can be done straight forwardly. However, increasing $k$ or $d$ seems to be very difficult.

Since it is possible that for large values of $n$ and $d$, there are no monochromatic graphs, we introduce a measure of monochromaticness on bi-colored graphs as follows. 

\begin{definition}[monochromatic Fidelity] Let $N$ be
$$N=\sum_{c} \left|w(c)\right|^2,$$
let $\mathcal{C}^{mono}$ be the set of all monochromatic IVC of $G$, and $d$ be the number of different colors of $G$.
The \textit{monochromatic fidelity} is defined as
$$F^{mono} :=\frac{1}{d} \frac{1}{N}\left|\sum_{c \in \mathcal{C}^{mono}}w(c)\right|^2.$$
\end{definition}
As an example, we can calculate the monochromatic fidelity of the graph in Figure \ref{fig:fig0MultiGraphVertex6}. It has $d$=3 monochromatic inherited vertex colorings and $N=6$. Then we find that $F^{mono}=\frac{3}{6}=0.5$. Furthermore, all monochromatic graphs reach the maximum of $F^{mono}=1$. 

\begin{Question}{approximative monochromatic graph}{} For every value of $n$ and $d$, which bi-colored graphs $G$ with $n$ vertices and $d$ different colors maximizes the monochromatic fidelity $F^{mono}$?\label{Question3}
\end{Question}

Even if one has access to $n-k$ heralding particles, it is possible that there are no $k$-monochromatic graphs with $d$ different colors, therefore we can define a fidelity as follows. 
   
\begin{definition}[$k$-monochromatic Fidelity] Let $N$ be
$$N=\sum_{c} \left|w(c)\right|^2,$$
let $\mathcal{C}^{k-mono}$ be the set of all $k$-monochromatic IVC of $G$, and $d$ be the number of different colors of $G$. The \textit{k-monochromatic fidelity} is defined as
$$F^{k-mono} :=\frac{1}{d} \frac{1}{N}\left|\sum_{c \in \mathcal{C}^{k-mono}}w(c)\right|^2.$$
\end{definition}    
   
For $k$-monochromatic states, the fidelity is $F^{k-mono}=1$. Naturally, we can ask what graph is closest to monochromatic.   

\begin{Question}{approximative $k$-monochromatic graph}{} For every value of $n$, $d$ and $k$, which bi-colored graphs $G$ with $n$ vertices and $d$ different colors minimizes the $k$-monochromatic fidelity $F^{k-mono}$?
\label{Question4}
\end{Question}

Until now, we considered only monochromatic colorings, as they correspond to an important class of quantum states. However, in general we are interested in the total capability of photonic quantum experiments to create quantum states. For that, we generalise our questions such that we cover every pure quantum state.

\begin{Question}{general inherited vertex colorings}{} 
Let $\mathcal{C}_p=\{C_i\}_{i=1}^{t}$ be a set of (prescribed) different colorings of $n$ vertices and $\mathcal{W}_p=\{w_i\}_{i=1}^{t}$ be the set of (prescribed) weights. For every $\mathcal{C}_p$ and $\mathcal{W}_p$, is there a bi-colored graph $G$ on the same $n$ vertices as the colorings in $\mathcal{C}_p$ so that for each $i$, $w(C_i)=w_i$, and every coloring not in $\mathcal{C}_p$ cancels out? 
\label{Question5}
\end{Question}

\begin{figure}[!ht]
\centering
\includegraphics[width=0.7\textwidth]{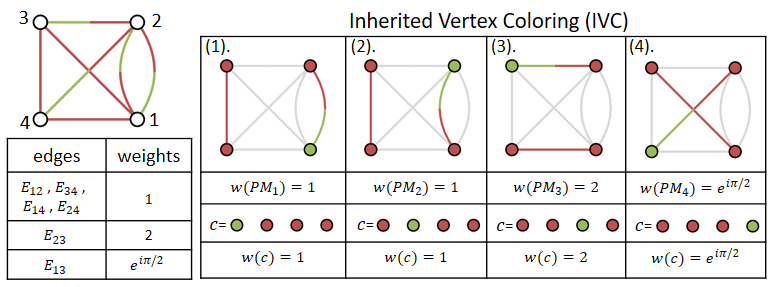}
\caption{This multi-edge graph answers the Question 5 for a given $\mathcal{C}_p$ and $\mathcal{W}_p$.}
\label{fig:fig5MultipleEdge}
\end{figure}

A particularly interesting special case of this question is the case where $\mathcal{C}_p$ is restricted to contain only $d=2$ colors. As an example, we consider the set of colorings $\mathcal{C}_p=\left((g,r,r,r),(r,g,r,r),(r,r,g,r),(r,r,r,g)\right)$ and weights $\mathcal{W}_p=(1,1,2,i)$. Is there a graph which is affirmative to Question 5 with these colorings and weights? We answer this question affirmatively, and show the solution in Figure \ref{fig:fig5MultipleEdge}.  

Again, it might be the case that not every set of coloring and weight can be constructed, thus we define a fidelity that gives us a notion of distance between the target and the graph. 

\begin{definition}[general fidelity]
Let $\mathcal{C}_p=\{C_i\}_{i=1}^{t}$ be a set of (prescribed) different colorings (with up to $d$ different colors), and $\mathcal{W}_p=\{w_i\}_{i=1}^{t}$ be the set of (prescribed) weights, let $G$ be a bi-colored graph. Let $N_1$ and $N_2$ be
$$N_1=\sum_{i=1}^{t} \left| w_i \right|^2, \quad N_2=\sum_{\forall c} \left|w(c)\right|^2.$$
The \textit{general fidelity} is defined as
$$F^{general} :=\frac{1}{N_1 N_2} \left|\sum_{i=1}^t w_i\cdot w(C_i) \right|^2.$$
\end{definition} 

Now a natural and most general question can be stated as follows.

\begin{Question}{approximative general graph}{}
For every $\mathcal{C}_p$ and $\mathcal{W}_p$, which bi-colored graphs $G$ with $n$ vertices minimizes the general fidelity $F^{general}$?
\label{Question6}
\end{Question} 

Question 6 contains Question 1-5 as special cases. Thus its resolution would resolve the question about the power of modern photonic quantum entanglement sources.

\section{Quantum Mechanical Formulation}\label{sec:QM}
All of the concepts, questions and partial results in this paper can directly be translated into the language of quantum mechanics \cite{krenn2017quantum,gu2018quantum,gu2018quantum3}.

\begin{description}
\setlength\itemsep{1em}
\item [Undirected Graphs] correspond to quantum optical experiments, using probabilistic photon-pair sources and linear optics.\par
\item [Vertices] correspond to single photon detectors in the output of some photon path.\par
\item [Edges] correspond to photon pairs that emerge from two photon paths.\par
\item [Edge weights] correspond to the amplitude of the corresponding photon pair.\par
\item [Edge colors] correspond to the mode number of the two photons in the path defined by the vertices at the endpoint of the edge. They can be bi-colored, as the two photons can have different mode numbers. A monochromatic edge corresponds to a photon pair with the same mode number.
\item [Perfect matchings] correspond to a multi-photon event where each single photon detector detects a photon. The coherent sum of all perfect matchings leads to the quantum state (conditioning on the click of each detector). Not every perfect matching necessarily leads to an unique term in the quantum state. Different perfect matchings can lead to the same inherited vertex colorings, thus coherently sum up and constructively or destructively interfere. 
\item [Inherited vertex colorings] corresponds to multi-photonic terms with different mode numbers in the quantum state. Terms with different IVCs are orthogonal.
\item [Weights of vertex colorings] $w(c)$ correspond to the amplitude of terms with mode numbers described by the inherited vertex colorings. More than one perfect matching can lead to the same inherited vertex colorings. As these terms can have opposite weights, it could be that the weight of an inherited vertex coloring is zero even though there are several perfect matchings leading to that coloring with nonzero weights.
\item [Monochromatic vertex colorings] lead to terms where every photon carries the same mode number. A graph with only monochromatic vertex colorings (with $d$ different colors) corresponds to $d$-dimensional Greenberger-Horne-Zeilinger state. These states are of significant importance in quantum physics.
\item [Question 1] asks which high-dimensional Greenberger-Horne-Zeilinger states can be created if general amplitudes $w_i \in \mathbb{C}$ can be used, but without trigger photons.
\item [Monochromatic Graph] corresponds to a high-dimensional multi-photonic Greenberger-Horne-Zeilinger state.
\item [Bogdanov's Lemma] states that Greenberger-Horne-Zeilinger states can be created only with $d=3$ dimensions with $n=4$ photons, or $d=2$ dimensions for arbitrary even number of $n$ photons, if all amplitudes are real valued (i.e. no destructive interference happens) and no additional trigger photons are used \cite{bogdanov267013}.
\item [Figure 2] corresponds to a 4-photon 3-dimensional Greenberger-Horne-Zeilinger state.
\item [$k$-monochromatic colorings] correspond to quantum states where the first $k$ photons have the same mode number, and the remaining $(n-k)$ photons have mode number zero (we can define $red$ to be an arbitrary mode number). The $(n-k)$ red vertices can be used as trigger photons that herald an $k$-photon state where every photon has the same mode number.  
\item [Question 2] asks which high-dimensional Greenberger-Horne-Zeilinger states can be created if general amplitudes $w_i \in \mathbb{C}$ can be used, and ($n-k$) trigger photons can be used.
\item [Erhard graph] is the only known example which corresponds to a quantum state that goes beyond Bogdanov's limit -- it can produce a 6-photon 3-dimensional entangled GHZ state. Four heralding photons and complex weights are used to cancel out all non-monochromatic colorings. It is created using two copies of the graph in Figure 2, which are merged using a quantum technique discovered by Manuel Erhard.
\item [Monochromatic fidelity] stands for a quantum fidelity to a high-dimensional $n$-particle GHZ state.
\item [Question 3] asks for every $d$-dimensional and $n$-particle state, what is the state that comes closest to the GHZ state, allowing only linear optics and probabilistic pair sources. 
\item [$k$-monochromatic fidelity] stands for a quantum fidelity to a high-dimensional $k$-particle GHZ state, using ($n-k$) trigger photons.
\item [Question 4] asks for every $d$-dimensional and $k$-particle state with ($n-k$) triggers, what is the state that comes closest to the GHZ state, allowing only linear optics, probabilistic pair sources and heralding photons. 
\item [Question 5] asks in general, which high-dimensional multipartite pure quantum states can be created using these techniques?  
\item [Figure 4] is an example to produce a 4-particle W state.  
\item [General fidelity] corresponds to a fidelity between a prescribed quantum state, and a quantum state that originates from a bi-colored graph.
\item [Question 6] asks for an arbitrary pure quantum state, with which fidelity can it maximally be created?
\end{description}

\section{Conclusion}\label{sec:conclusion}
Every progress in any of these purely graph theoretical questions can be immediately translated to new understandings in quantum physics. Apart from the intrinsic beauty of answering purely mathematical questions, we hope that the link to natural science gives additional motivation for having a deeper look on the questions raised above.

\addtocontents{toc}{\SkipTocEntry}
\section*{Acknowledgements}
The authors thank Manuel Erhard, Anton Zeilinger, Tomislav Do\v{s}li\'{c} and Roland Bacher for useful discussions and comments on the manuscript. M.K. acknowledges support from by the Austrian Academy of Sciences (\"OAW), by the Austrian Science Fund (FWF) with SFB F40 (FOQUS). X.G. acknowledges support from the National Natural Science Foundation of China (No.61771236) and its Major Program (No. 11690030, 11690032), the National Key Research and Development Program of China (2017YFA0303700), and from a Scholarship from the China Scholarship Council (CSC). D.S. acknowledges support by the National Research, Development and Innovation Office NKFIH, No. K-120706, No. KH-130371 and No. KH-126853.

%\begin{thebibliography}{10}
\bibliographystyle{unsrt}

\bibliography{refs}

\end{document}